Dmytro Taranovsky
July 17, 2017

# Finitistic Properties of High Complexity

**Abstract:** We use fast-growing finite and infinite sequences of natural numbers and more complicated constructs to define models of hypercomputation and interpret non-arithmetic predicates, with the strongest extensions reaching full second order arithmetical truth and beyond. Since the predicates are interpreted using properties of certain natural finite structures, they are arguably finitistic.

**Contents:**

## 1 Introduction

Finitism refers to definitions and arguments that do not use an actual infinity, though they may involve a potential infinity. The potential infinity refers to the fact that the natural numbers are unlimited, and that there are arbitrarily large natural numbers. Hyperfinitism by contrast deals only with "small" natural numbers, such as those whose binary representation is feasible.

The central question about finitism is which properties are definable and which claims are demonstrable finitistically. At one extreme is the view that finitistic properties are the recursive ones, and that finitistic demonstration does not go much further beyond PRA. At the other extreme is the view that every definable property (including the truth predicate of set theory) of natural numbers can be defined finitistically, and that every true mathematical claim is finitistically demonstrable. On that view, infinite structures are useful as conceptual aids and for argument simplification, but can ultimately be replaced by finitistically definable finite analogues. Finitistic strong consistency proofs can involve either new finitistic properties or just new insights about the recursive properties. Most logicians hold a narrow view of finitistic definability because of absence of evidence to the contrary. Plausibly finitistic ways to define properties of high complexity or to prove strong consistency statements were simply not known before this paper (or at least not as generally known or developed). One compromise view is that philosophically, there is a continuum of notions between being fully finitistic and using infinity, and our results help to elucidate that continuum. And one formalist (or under Platonism, metaphoric) view is that our notions are an alternative approach to infinity: Infinity can be expressed directly (using infinite sets) or indirectly as idealization of large numbers and fast growing sequences and other constructs (as this paper develops). In any case, whether or not the definitions here are finitistic, they enrich our foundational understanding by interpreting statements involving infinity as properties of certain natural finite structures.

Hypercomputation refers to idealized computers computing non-recursive functions. Models of hypercomputation can be divided into two types. One uses oracles, either directly or through some auxiliary mechanism. The other uses an infinite amount of computation in finite time. In the models we propose (which are technically oracle based), finite computations correspond to infinite computations by using fast-growing sequences to convert unbounded searches into bounded ones.

In a way, this paper complements my paper [1], especially the last section, and its follow-up [2]. These papers use largeness to extend the language — of number theory in this paper, and of set theory in the other papers — though here the extensions are definable in preexisting mathematics.

For each notion here, the strategy is to intuitively define the notion using just finite numbers, give a formal definition using infinite sets, argue that the intuitive and formal notions agree, and prove the complexity of the formal notion. Precise notions will be intuitively defined by defining a vague notion and then defining the precise notion such that the definition does not depend on how vagueness is resolved. In the infinitary formal treatment, possibilities for the vague notion will form a directed system such that every sufficiently large choice works correctly.

The second section describes computation using a single sufficiently fast-growing infinite sequence. The third section discusses computation using multiple sequences. The fourth section presents similar results for finite sequences. The fifth section describes an arguably finitistic way to compute $\Sigma^1_2$ truth. The sixth section uses estimators to extend the results to full second order arithmetic, and the final section goes beyond second order arithmetic using transfinite levels of estimators (7.1), pattern enumerators (7.2), and countable infinity (7.3). We also uncover connections for different levels of notions, including connections with games and determinacy at the levels of open determinacy (for a single sequence), determinacy for the $\Sigma^0_2$ difference hierarchy (for multiple sequences), projective determinacy (using estimators), and beyond (using pattern enumerators).

## 2 A Model of Hypercomputation

A powerful model of hypercomputation rests on a surprisingly mild physical assumption:

> There is a type of event that occurs an unlimited number of times, but with frequency of occurrence decreasing sufficiently fast.

**Definition 2.1:** A language $L$ is recognized by a Turing machine with *the fast-growing sequence oracle* iff there is a total function $A$ such that for every $B$ with $B(n) \geq A(n)$ (for every natural number $n$), the machine using $B$ as an oracle halts on input $S$ iff $S$ is in $L$.

A variation on the definition would be to require $A(n)$ to be sufficiently large relative to $(A(0), A(1), \ldots, A(n-1))$, but the expressive power is the same as we can get the stronger version using $(A(0), A(A(0)), A(A(A(0))), \dots)$. In this section, we will use "for every sufficiently fast-growing $B$, $\varphi(B)$" to mean $\exists A \forall B (\forall n B(n) \geq A(n))\, \varphi(B)$, but the *strict version* of "sufficiently fast-growing" $B$ requires
$B(0) \geq A(0) \wedge \forall n B(n+1) \geq A(B(n))$. Also, for us, it is insufficient for $B$ to just be fast-growing asymptotically since the Turing machine must be chosen independent of $B$.

It is not known whether such machines are physically constructible. Currently, the biggest problem is surviving for the extremely long time required, rather than the apparent lack of the required physical events. In any case, studying the machines improves our understanding of the recursion theory. Below, we prove that a language is recognized by such a hypercomputer iff it is $\Pi^1_1$. Thus, a language is decidable in the model iff it is hyperarithmetic. Note that $\Pi^1_1$ rather than $\Sigma^1_1$ is the correct analogue of recursively enumerable.

**Theorem 2.2:** Recognizability for Turing machines with the fast-growing sequence oracle equals $\Pi^1_1$.
**Proof:** For a sufficiently fast-growing sequence $A$, a recursive relation has an infinite descending path through $n$ iff it has an infinite descending path through $n$ and then through a natural number less than $A(\max(n,m))$ where $m$ is the length of the given recursive definition of the relation. By König's lemma, if the relation is well-founded, then the tree is finite, and thus the machine searching it will eventually discover that the relation is well-founded.
For the converse, for each machine and input in the theorem, the machine halts iff for every $A: \mathbb{N} \to \mathbb{N}$ there is a halting computation such that the answers the machine receives are at least as large as the answers given by $A$.

If we do not require a Turing machine to act uniformly on all sufficiently fast-growing $g$, we have the following.

**Theorem 2.3:**
*(a)* The set of Turing machines that halt for some sufficiently fast-growing $g$ (that is $g(n)$ is sufficiently large relative to $n$ and the Turing machine) is $\Pi^1_1$-complete.
*(b)* The set of Turing machines that halt for every monotonic sufficiently fast-growing $g$ is $\Pi^1_1$-complete.
*(c)* The set of Turing machines that halt for every sufficiently fast-growing $g$ equals (under many-one reduction) the game quantifier for $\Sigma^0_2$-games on integers. This also applies to $\Pi^1_1$ formulas in place of Turing machines.
*Notes:*
\* The theorem (and the proof) relativizes to a real parameter.
\* As a corollary, (b) also holds for the class of $g$ with $g(0)$ sufficiently large and $g(n+1)$ sufficiently large relative to $g(n)$.
\* Game quantifier is the set of codes of games (given a coding) for which the first player has a winning strategy.
\* Lack of monotonicity makes the condition in (c) less natural. The bound in (c) is far beyond a finite iteration of the hyperjump.
**Proof:**
*(a)* Follows from the proof of Theorem 2.2.
*(b)* Given a Turing machine $T$, consider the following game:
Player II: $N$ (a number)

repeat:
   Player I: next $N$ values of a sequence $f$
   Player II: next $N$ values [alternatively, one value] of a monotonic sequence $g$ with each $g(n) \geq f(n)$; new value of $N > 0$.
First player wins iff $T$ halts using $g$.
This is an open game on integers, so the game quantifier (the set of $T$ for which player I wins) is $\Pi^1_1$. If $T$ halts for every sufficiently fast-growing monotonic $g$, then player I wins (by playing a sufficiently fast-growing $f$). Conversely, if for some sufficiently fast-growing monotonic $g$, $T$ does not halt, then player II wins because for every strategy for player I, player II can make $g$ into any monotonic function that is sufficiently fast-growing relative to the player I strategy. (Proof sketch: Given a well-behaved bound $h$ derived from player I strategy, and a desired $g(n) > h^{3n}(0)$, choose the least $N$ to get to a gap $g(n+1) > h(h(g(n)))$. This works since $(g(n), N)$ is sufficiently small relative to $g(n+1)$.)

*(c)* To compute a $\Sigma^0_2$ ($\exists m \forall n\, \varphi(m,n)$)-game quantifier, modify the game by requiring each move to be below g(2*number_coding_the_play_so_far), and modify the winning condition to $\exists m \forall n < g(2m+1)\, \varphi(m, n, \text{play})$, and halt if the first player has a winning strategy. Here we use $g$ being sufficiently large for even numbers, the possibility of $g$ for odd numbers being sufficiently large relative to $g$ for even numbers, determinacy, and König's lemma.
Conversely, a Turing machine $T$ halts for every sufficiently fast-growing $g$ iff the first player wins the following game: At move $n$, the first player plays $f(n)$ (thus defining $f$), while the second player may pass or play $g(m) \geq f(m)$ where $m$ is the least such that $g(m)$ has not yet been determined. The first player wins if $T$ halts using $g$, or some $g(n)$ is undefined. Given a first player strategy, the second player can make $g$ into an arbitrary sufficiently fast-growing function by delaying playing $g(n)$ until it is not too large relative to the move number.
Also, given the universal quantification, a $\Pi^1_1$ formula can be coded by a Turing machine that halts once the proposed function is inconsistent with being a counterexample.

For readers familiar with forcing, another interesting quantifier is "for every generic sufficiently fast-growing sequence" which is equivalent to "for every Hechler real (as a sequence of integers) with a sufficiently long initial segment removed". For non-absolute properties, the two versions are (1) $V$-generic, and (2) (in line with this paper) generic for a sufficiently complete countable model (using Template Definition 4.2) but existing in $V$. Hechler reals satisfy every open property (from the list of properties in the ground model) that is dense on sufficiently fast-growing sequences. To the extent determinacy holds, all such reals have the same properties for properties (in the ground model) that are not altered by removing a finite initial segment of the sequence. The sequences are not monotonic (in a sense anything that can happen happens), but we can define variants for monotonic, strictly monotonic, and the strict version (as used beneath Definition 2.1).

**Theorem 2.4:** A predicate is arithmetically definable from every sufficiently fast-growing sequence iff it is recursive in a finite hyperjump of 0.
**Proof:**
   The computation of hyperjump relativizes, with each extra quantifier capturing another hyperjump.
   For the converse, we show that for every arithmetic $\varphi$, $\{n : \forall G\, \varphi(n, G)\}$ is recursive in a finite hyperjump of 0, where $G$ ranges over Hechler reals. (Furthermore, non-uniformity in φ would not help, as there would still be a forcing condition that forces uniform definability.) Thus $G$ is (but for an arbitrary initial segment) a generic sufficiently fast-growing sequence. The conditions for Hechler forcing are (in one formulation) $(m, f)$ ($m \in \omega$ and $f \in \omega^\omega$) where $(n, g)$ is as at least as strong as $(m, f)$, denoted $(m, f) \geq (n, g)$ iff $m \leq n$ and $\forall i < m\, f(i) = g(i)$ and $\forall i\, f(i) \leq g(i)$. Using basic properties of forcing, a formula holds for every $G$ if every condition forces it. A condition $p$ forces $\forall n \in \omega\, \psi(n, G)$ iff it forces $\psi(n, G)$ for every $n$. A condition $p$ forces $\neg \psi$ ($\psi$ a statement that may use $G$) iff no stronger condition forces $\psi$.

To get the required complexity, let us say that fast-growing $f : m \to \omega$ forces a statement iff the statement is forced by $(m, f')$ for every sufficiently fast-growing $f' : \omega \to \omega$ extending $f$. (This implies, for example, $f'(m) \gg \max(f)$. Also, $f$ itself is not required to be fast-growing.) If we care only about a set of statements whose number is bounded relative to $m$, there are particular conditions equivalent to fast-growing $f : m \to \omega$, and furthermore (letting $f$ and $m$ vary) the set of such conditions is dense. Combining all together, a fast-growing $f : m \to \omega$ forces $\forall n\, \psi(n)$ iff there is a rate of growth $f'$ such that for every $n$ every fast-growing $g$ forces $\psi(n)$ provided that $g$ extends $f$ and $g(i) \geq f'(i)$. A fast-growing $f : m \to \omega$ forces $\neg \psi$ iff there is a rate of growth $f'$ such that no fast-growing $g$ forces $\psi$ provided that $g$ extends $f$ and $g(i) \geq f'(i)$. By induction on the complexity of the formula, this has the required complexity, the key point being that $f'$ plays no role in checking whether fast-growing $g$ (which is finite) forces a statement.

**Computational Complexity**

For each machine, input, and the set of possible $A$, we may define its computation tree, where we start at the root, each vertex is a query to $A$, and each branch is an answer to the query, and we may define the ordinal given machine and input as the ordinal depth of the computation tree for the set of sufficiently fast-growing $A$ (where 'sufficiently fast-growing' may depend on the machine and input). Intuitively, the ordinal says how long in the worst case the computation will take relative to the rate of growth of $A$. From there, the ordinal for a machine is the supremum of ordinals across inputs (assuming the machine always halts for sufficiently fast-growing $A$), and the ordinal for a problem is smallest ordinal for a machine (whose output is independent of the choice of a sufficiently fast-growing $A$) that solves the problem.

**Theorem 2.5:** For each recursive ordinal $\alpha$, the ordinal for a problem is at most $\alpha$ iff the problem is recursive in $0^{(\alpha)}$ (αth Turing jump of 0).

**Proof:** To compute $0^{(\alpha+1)}$, given an input $n$, compute $0^{(\alpha)}$ up to $A(n)$, and then check whether machine $n$ halts on $0^{(\alpha)}$ before getting out of bounds. The computation depth is $\alpha + 1$, as required (for limit $\alpha$, $0^{(\alpha)}$ can be computed as a sequence of lower Turing jumps).

Conversely, given a computation tree of depth $\alpha$ for sufficiently fast-growing sequences, start with accepting halting positions (without checking rate of growth), and work backwards, and at each stage, mark positions such that for every sufficiently large answer, the resulting position was marked. For $\alpha = \beta + n$, note that each level adds a $\Delta_2^0$-quantifier, so $n$ levels merge into $\Delta_{n+1}^0$ as required. (Each quantifier is $\Delta_2^0$ since assuming the sequence so far was sufficiently fast-growing, and given that we measure at maximum depth, every sufficiently large choice for the next element will give the same acceptance value.)

Informally, we can speak of computational complexity on a much more fine-grained level as follows: A problem can be solved in time $T(n)$ if there is a Turing machine and a natural set $A$ such that the machine using $A$ solves the problem in time $T(n)$ and that if $\forall n\, B(n) \geq A(n)$, the machine using $B$ correctly solves the problem. Thus, recursive languages have recursive complexity, arithmetic languages have arithmetic complexity, and so on. The converse might be false, and it is unclear what the complexity of the inverse Busy Beaver function should be. Unfortunately, we do not know how to define complexity formally because using unnatural $A$ we can code solutions to recursive problems as follows: If instance $k$ is true, set $A(2k)$ to a certain value, and otherwise make $A(2k)$ larger, and similarly with $k$ being false and $A(2k+1)$. One mitigation (which still breaks down at high complexity) would be to make $A$ sufficiently sparse, for example by requiring that $A$ is coded through a binary tape with sequential access such that every '1' is at a position $2^n$ (for some n) and that if a '1' is position $2^n$, the next '1' cannot occur until position $2^{n^2}$. That way, for problems with easily verifiable solutions there is at most a polylogarithmic speed up (relative to the running time) as we can try all possible initial oracle tape segments.

## 3 Strengthenings to Multiple Sequences

Still higher complexity is reachable by computers that use a sufficiently fast-growing function $A$ and a function $B$ that grows sufficiently fast relative to $A$. We can extend the model to allow more sequences, each of which growing sufficiently fast relative to the previous ones. We can even (not part of definition 3.1) allow transfinitely many sequences through the use of diagonalization.

**Definition 3.1:** A language L is recognizable by a Turing machine with the fast-growing level $k$ sequence oracle iff there is a Turing a machine T and a function $f : \mathbb{N}^\mathbb{N} \to \mathbb{N}^\mathbb{N}$ with f(A)(n) depending only on a finite (but dependent on $n$) segment of $A$ such that for every $A_1, \ldots, A_k$ with $\forall i < k\, \forall n\, A_{i+1}(n) > f(A_i)(n)$ ($A_0$ is 0 here), T halts on input S with oracles $A_1, \ldots, A_k$ iff S is in L. The *nonlocal version* allows $f(A)(n)$ to depend on the entirety of $A$.

*Notes:*
* We use the term "level $k$" (or $k$-level) because for a number $n$ in the range, the maximum $i$ such that $n \in A_i$ is from 1 to $k$.
* As before, a strict version of fast-growing level $k$ sequences (that does not affect Definition 3.1) is to require $A_{i+1}(n+1) > f(A_i)(A_i(n))$ ($A_i(-1)$ being 0 here); one could also require $A_{i+1} \subset A_i$ (with each $A$ viewed as a set). The expressive power remains the same, as one could simulate Turing machine that assume these more restrictive versions, including a Turing machine that assumes the growth rates are sufficiently high relative to the input size. This applies to both local and nonlocal versions.
* The definition can be extended to transfinitely many levels, with $A_i(n) > f((A_j : j < i))(n)$.
* We conjecture that local and nonlocal versions have the same expressive power (see

Theorem 3.2, and Conjecture 3.3). The local version is more metaphysically parsimonious, and has the same expressive power as if we treated a fast-growing infinite level $k$ sequence as an infinite sequence all of whose sufficiently long initial segments are (in the sense of Definition 4.3(c)) sufficiently long fast-growing level $k$ sequences. However, quantification over sequences depends on the version. For example, whether a Turing machine halts for some sufficiently fast-growing level 2 sequence (nonlocal version) is $\Sigma^1_2$ complete.

**Theorem 3.2:** A language is recognizable by a Turing machine with the fast-growing level 2 sequence oracle iff it is many-one reducible to the $\Sigma^0_2$ game quantifier (equivalently, to the set of winning positions in a $\Sigma^0_2$ game on integers). This also applies to the nonlocal version.
**Proof:** If a player has a winning strategy, then he has a winning strategy using moves less than $A_1$(position), where position includes the position in the game and a canonical description of the rules. Given this constraint on the moves, the second player has a winning strategy iff he has a strategy that enters then $n$th open set before time $A_2(n+\text{position})$. Given these constraints, the game becomes an open one with finite number of moves per turn, and hence recognizable.
For the converse, if a Turing machine $T$ behaves uniformly (in either halting or not halting) using a fast-growing level $k$ sequence oracle, then $T$ halts on input $s$ iff the first player wins the following game:
Player I: $A_1(n)$
Player II: $B_1(n) \geq A_1(n)$, optionally $(m, A_2(m))$
Player I wins iff $T$ halts on $s$ using $B_1$ and some $B_2$ with $\forall m B_2(m) \geq B_1(m)$, or for some $m$, $A_2(m)$ is undefined or inconsistent.
Since the winning condition depends only a finite segment of $B_2$, it is $\Sigma^0_2$. Player I can make $B_1, B_2$ (for some choice of $B_2$) a sufficiently fast-growing level 2 function. Conversely, given a player I strategy, player II can simultaneously make $B_1$ an arbitrary sufficiently fast-growing function (strict version) and pick an arbitrary $A_2$, and thus able to guarantee that $B_1, B_2$ is a sufficiently fast-growing level 2 function.

**Conjecture 3.3:** We conjecture the following to have equal complexity ($n > 0$):

1. Recognizability using sequences of level $n+1$.
2. Recognizability using sequences of level $n+1$ (nonlocal version).
3. The game quantifier for level $n$ of the difference hierarchy of (lightface) $\Delta^0_3$ sets (level 1 is $\Sigma^0_2$, level two is $\Sigma^0_2 \setminus \Sigma^0_2$, level 3 is $\Sigma^0_2 \setminus (\Sigma^0_2 \setminus \Sigma^0_2)$, and so on).
4. Monotone inductive definitions for co-recognizable sets (relative to the current real number) for sequences of level $n$ (for $\Pi^0_1$ if $n = 0$).

Moreover, we conjecture that allowing number quantification in the model (i.e. allowing arbitrary arithmetical formulas in place of Turing machines) corresponds to boldface (in terms of $\Sigma^1_1$ strength of the corresponding theory) or finitely iterated versions of these properties.

Furthermore, using the results in [4], we appear to have the following equivalence:
*(a)* $x$ is expressible using (3) (in Conjecture 3.3) for some $n$.
*(b)* $x$ is in the least $L_\kappa$ satisfying $\forall n \exists \kappa_1, \ldots, \kappa_n \, \kappa_1 < \ldots < \kappa_n < \kappa \wedge L_{\kappa_1} \prec_1 \ldots \prec_1 L_{\kappa_n}$
(these sets for each numeral $n$ and relativized to each real $y$ exactly capture the $\Pi^1_3$

strength of $\Pi_2^1 - \mathrm{CA}_0$).
*(c)* $x$ is recursively reducible to the set of trees accepted by a tree automaton (equivalently by a formula in S2S) (the reduction $f$ is such that $x(i) = 1$ iff $f(i)$ is a (possibly infinite) tree that is accepted).

**Towards Proving the Conjecture**

**Theorem 3.4:** The set of Turing machines that halt for every (alternatively, some) sufficiently fast-growing growing level $k > 1$ sequence (strict version) is many-one reducible to the game quantifier for level $k - 1$ of the $\Sigma_2^0$ difference hierarchy (as in Conjecture 3.3(3)).
*Note:* Using theorem 3.5, the converse also holds.
**Proof:**
We prove that the machine halts for every such sequence iff player I wins the following game:
player I: $A_1(n)$, and $(i, m, A_i(m))$ for each undefined $A_i(m) \leq A_1(n)$ ($1 \leq i \leq k$).
player II: $B_1(n)$, and $(i, m, B_i(m))$ for each undefined $B_i(m) \leq B_1(n)$ ($1 \leq i \leq k$); $B_i(m) \geq A_i(m)$ (and cannot be played until $A_i(m)$ is defined).
Player I wins iff: (1) some $B_i(m)$ is undefined and for all $(i', m') \leq (i, m)$ (lexicographically) $A_{i'}(m')$ is defined, or (2) $T$ halts using $B_1, B_2, \ldots, B_k$ (querying an undefined value waits without halting until it is defined).
Given an $f$ (as in Definition 3.1), player I can force $B$ to be consistent with $f$. For the converse, fix a strategy for player I, and without loss of generality, if $T$ halts, set all undefined $A_i(m)$ to some values. If for some choice of player II moves, player I fails to produce some $A_i(m)$ with $B_{i'}(m')$ defined for all lexicographically lower $(i', m')$, player I loses, so assume otherwise. With some basic assumptions, for a fixed complete $B_i$, there is a finitely branching tree of possibilities (with each vertex corresponding to an entry in $B_i$), and set $A'_{i+1}(B_i)(n)$ to be the largest possible value of $A_1$ at which player II can learn the next value of $A_{i+1}$ if the current value of $B_{i+1}$ is defined and is $\leq n$. Using $A'$, one can choose $f$ so that player II can get an arbitrary $B_1, B_2, \ldots, B_k$ consistent with $f$, and thus win unless the machine halts for every such $f$.
For the version where the machine halts for some such sequence, the game is similar but with the second player winning if the machine halts, and the game quantifier for the second player has the required complexity.

**Theorem 3.5:** A language is recognizable by a Turing machine with the fast-growing level $k > 1$ sequence oracle iff it is many-one reducible to the game quantifier for level $k - 1$ of the $\Sigma_2^0$ difference hierarchy.
**Proof:**
   Theorem 3.4 proves (a strengthening of) one direction. In the other direction, games on integers with the payoff in level $k - 1$ of the $\Sigma_2^0$ difference hierarchy can be coded by games on $\{0,1\}$ with the payoff in level $k$. Furthermore, the payoff can be converted to the following form: First player wins iff the highest event priority seen infinitely often is even for odd $k$ and odd for even $k$ (where event level 0 occurs on every step). Events will be observable (computable), and a higher priority event is also counted as an event of all lower priorities. Fix such a game, and for each position, a winning strategy for the right player. We show that the game is equivalent to the following modification: If after time $n$, a priority $i$ event is not seen until $A_i(n)$ then the player (who wants level $i$ to be seen infinitely often) loses. (The loss happens at time $A_i(n)$, and if the losing condition is met for multiple $i$ (with $A_i(n) = A_{i'}(n')$), the lowest one controls.) This is an open game on {0,1} and hence recognizable.

Set $f_1(n)$ (in any specific way) such that if the player (who wants an odd level to be seen infinitely often) wins starting at some position at time $n$, then using the strategy we fixed, and regardless of the opponent's moves, the next event (or a victory) will be seen before time $f_1(n)$. Since $A_1(n) \geq f_1(n)$, the game and strategies are unaffected by the modification with $A_1$. Now compute $f_2(A_1)(n)$ using the same strategies, events of priority 2, and the modified winning condition. (Technically, we have a set of games, one per position; if the position is at time $n$, the constraint only uses $A_1(n)$ or higher.) Note that $f_2$ can be chosen to only depend locally from $A_1$, and furthermore, we can detect a bad $A_1(n) < f_1(n)$ and ensure that $f_2(A_1)$ remains defined. Continuing further, we get all $f_i$ (using bounded maximization to depend only on $A_i$ and not $A_1, \ldots, A_{i-1}$), and combine them into $f$ as in Definition 3.1.

It appears likely that for transfinite levels, recognizability for a level α sequence corresponds to the game quantifier for level α of the difference hierarchy. In these games, events can have α different levels, and if the lowest level event not seen infinitely often is a limit, the second player wins iff α is even (limit ordinals are even here).

Relation with monotonic induction:
* The game quantifier can be computed using monotonic induction: Use the previous quantifier to find positions where player one can win while preventing the game from getting into a state labeled with $k$ (where player two wins if $k$ is seen infinitely often), and mark these positions as player one wins, and repeat with the modified game, continuing using monotonic induction until no new positions are added. The other positions are player two wins. (Side note: A winning strategy for player one would be to play the auxiliary game for the step where the current position was marked, and each time a previously marked position is reached, player one wins the auxiliary game and switches to the auxiliary game that led to the position being marked. For player two, the auxiliary game (this is not the auxiliary game for player one) is keeping to unmarked positions, with a win if a state labeled with $k$ is reached. For both players, if the auxiliary game does not halt at a finite step, the winning condition is the same as in the original game.)
* To simulate the multiple sequence oracle (conjectured), we can build the computation tree where each vertex corresponds to a query to $A_k$. Work backwards from halting configurations, and using the previous quantifier for each step, use monotonic induction to reach the initial state if it is halting.

Relation with nonlocal version and with arithmetic definability: We conjecture that by using a $k$-fold iteration of Hechler forcing, we get fast-growing level $k$ sequences that give minimal uniform definability.

Going much further (than level $k$ sequences) but remaining inside a fragment of $\Delta^1_2$, an interesting nonfinitary notion (explored in literature) is that of infinite time Turing machines. Such a machine is like an ordinary Turing machine, except that it can run for an unlimited ordinal length, and at limit steps, the state of each cell (and the machine state) becomes lim inf of the previous states.

# 4 Finite Structures

The infinite structures have finite analogues. A key concept is that of a sufficiently fast-growing sufficiently long finite sequence. Sufficiently long is with respect to the rate of growth, and is preserved under extending the sequence. In computation, along with input, we will be given such a sequence, with the first element sufficiently large relative to the input size. Vagueness will be avoided because there is no error in interpreting "sufficiently" too strictly. What we are defining are not the sequences themselves, but the results of the computations — for computations that agree for all appropriate sequences (and thus will still agree if "sufficiently" is interpreted too strictly). The appropriate sequences can then be formally defined as the ones giving the correct results. Also, by using infinite sets, we can quantify over all possible definitions (like we did for the fast-growing sequence oracle) and thus formally define the results without using "sufficiently".

Let us now describe the construction in detail. In the definitions, 'vague' indicates dependency on a notion of 'sufficient'.

**Vague Definition 4.1:**
*(a)* A *sufficiently long fast-growing sequence* is a finite function $g : \{0, \ldots, N\} \to \mathbb{N}$ such that
  - $g(0)$ is sufficiently large (relative to the problem size),
  - $g(n+1)$ is sufficiently large relative to $g(n)$ ($n < N$), and
  - $g$ or an initial segment of $g$ is sufficiently long relative to its rate growth. The rate of growth is how large is $g(0)$, and how large is $g(n+1)$ relative to $g(n)$ (for $n+1$ in the range of the initial segment).
*(b)* A sufficiently long fast-growing sequence or its extension is any finite extension of a sequence in (a).
*(c)* (generalization of (a)) A sufficiently long fast-growing level $k+1$ sequence is a sufficiently long finite sequence of sufficiently long fast-growing level $k$ sequences such that once the sequences are expanded to numbers, each number is sufficiently large relative to the previous one. Here, a level 0 sequence is a number sufficiently large relative to the problem size. "Sufficiently long" may depend on the sequence (specifically, required length is relative to the rate of growth or its analog) but if a sequence is sufficiently long, then so are all of its extensions.
*Note:* One useful representation of a level $k$ sequence is a function $g : \{0, \ldots, n\} \to \{0..k+1\}$ with $g(x) \geq i$ iff $x$ is the last element (after expansion into a sequence of numbers) of an included level $i-1$ sequence (and $g(n) = k+1$; alternatively set $g(n) = k$ since with $g(n) = k$ extensions of the sequence are represented by extensions of $g$). For example, (((1,2,3),(4,5,6,7)),((10,11),(12,13,14))) would be represented as (0,1,1,2,1,1,1,3,0,0,1,2,1,1,4).

In one finitistic treatment, we assume that in certain contexts, 'sufficiently' is a primitive concept — it can be explained, partially axiomatized, and connected to other concepts, but not formally defined in terms of more primitive concepts. (On the other hand, one could argue that our uses of 'sufficiently' and other terms are a recharacterization of aspects of infinity.) Using that concept, given a particular Turing machine T, we can define a predicate $P$ as $P(n)$ holds iff $T$ accepts given input $n$ and a sufficiently large/complete/etc example (relative to $n$ and $T$) of the appropriate type, provided that $T$ is a machine that halts with the same (dependent on $n$) output for every $n$ and sufficient example. While without context, 'sufficient' is vague, the above usage is (at least arguably) precise because it does not depend on how the vagueness is resolved. One would then finitistically study the properties of $P$ and other such predicates, and develop axiom systems. However, with a little reasoning, many uses of "sufficiently" can also be formally defined using infinite sets, and we will use infinity as a shortcut and to unambiguously interpret the concepts in existing mathematics.

**Template Definition 4.2.**
*Assumption:* We are given a set of potential examples, and a set $S$ of possible notions of sufficiency: The notions must form a directed set under strictness, and each notion must be nonempty, or the definition fails. Alternatively, instead of $S$, we may be given $T$ that is equivalent to $S$ in that for every notion in $S$ there is at least as strict notion in $T$ and vice versa.
*(a)* A property holds for every sufficient example iff there is a notion $Q$ in $S$ such that the property holds for every example in $Q$.
*(b)* A predicate $P$ is computable using a sufficient example (assuming the examples are reducible to integers) iff there is a Turing machine $T$ such that for every $n$, $T$ outputs $P(n)$ given $n$ and a sufficient example (that is for every sufficient example; sufficiency depends on $n$) as input.
*(c)* Given a well-defined language, a sentence using a sufficiently good example as a parameter is well-defined iff there is a valid notion $Q$ such that all elements of $Q$ agree on whether the sentence is true, which is then the truth value of the sentence.

To apply the definition (which we will sometimes do implicitly), we need the set of possible notions of sufficiency. Here is one choice, noting that selection between equivalent (as used in 4.2) sets may be arbitrary.

**Definition 4.3:**
*(a)* A notion of sufficiently long fast-growing sequences is valid iff there is a function $f$ such that for every infinite sequence $g$ with $g(0) > f(0)$ and $g(n+1) > f(g(n))$, there is $m$ such that every initial segment of $g$ of length at least $m$ is sufficient.
*(b)* The extension for 4.1(b) is trivial. (A notion of sufficiently long fast-growing sequence or extension thereof is valid if it is the closure under finite extension of a valid notion in (a).)
*(c)* A notion of sufficiently long level $k+1$ (fast-growing) sequences is valid iff there is a function $f$ and a valid notion $P$ of sufficiently long level $k$ sequences such that for every infinite sequence $g$ of elements of $P$ with, when $g$ is expanded into a sequence of numbers, $g(0) > f(0)$ and $g(n+1) > f(g(n))$, there is $m$ such that every initial segment of $g$ of length at least $m$ is sufficient. For every $n$, $\{x : x \geq n\}$ is a valid notion of level 0 sequences.

**Theorem 4.4:**
*(a)* A predicate is computable using a sufficiently long fast-growing sequence iff it is $\Delta_1^1$ in the hyperjump of 0.
*(b)* A predicate is computable using a sufficiently long fast-growing sequence or extension thereof iff it is $\Delta_1^1$.
**Notes:**
* There should also be an analogous statement for 4.3(c).
* A finite sequence (as in (a)) is more expressive than the corresponding infinite sequence because it gives us a stopping point.
* For mass problems (where the correct answer is non-unique), the notion in 4.1(b)

can be used to compute $\Pi^1_1$ separation (with the choice of the correct answer depending the sequence). A weakening of 4.1(b) that does not allow that would be to store the sequence as a predicate so one cannot test whether one reached the end.
**Proof:**
*(a)* The algorithm above (Theorem 2.2) for testing well-foundness (with reaching the end of the sequence indicating an infinite path), will work correctly and compute the hyperjump of 0 up to a sufficiently large unspecified point in the sequence. In deciding a $\Delta^1_1$ in the hyperjump of 0 problem, use the partially computed hyperjump of 0 to check which of the two opposing sequences (for the sentence and its negation) terminates first: For inputs in the size limit, that will happen before that unspecified point.
Conversely, a machine that behaves uniformly halts iff for every growth rate $A$, there is a finite sequence $h$ growing faster than $A$ such that for every sufficiently-fast growing finite extension of $h$, the machine halts. Extensions of $h$ can be considered independently of $A$, so this is $\Pi^1_1$ in the hyperjump of 0.
*(b)* $\Delta^1_1$ predicates are clearly computable. For the converse, for every rate of growth $g$, there is a finite sequence $h$ with $h(n+1) > g(h(n))$ such that the output is given on every continuation of $h$ (and one similarly has a $\Sigma^1_1$ description).

Without uniformity, we have the following.

**Proposition 4.5:** For every $\Sigma^1_1$ $\varphi$, $\{n : \forall S \, \varphi(n, S)\}$ is $\Pi^1_1$ in the hyperjump of 0 where $\forall S$ quantifies over sufficiently long fast-growing sequences.
**Proof:** Consider the following game ($i$ is move number, $n$ is fixed):
First player: $A_1(i)$.
Second player: $A_2(i) \geq A_1(i)$.
At any point, either player can request to stop.
The first player wins iff once both players request to stop, $\varphi(A_2)$ holds, or the first player requests to stop, but the second player never does.
This game is a concatenation of an open game (before the first player asks to stop) and a closed game afterwards (also, φ can be tested using a closed game) without extra parameters, so the game quantifier is $\Pi^1_1$ in the hyperjump of 0. The first player can force the $A_2$ to be some sufficiently long fast-growing sequence, and relative to a strategy for the first player, the second player can force $A_2$ to be an arbitrary sufficiently long fast-growing sequence, so the first player wins iff $\forall S \, \varphi(n, S)$ holds above, as desired.

If we allow the computers to produce sufficiently long level $k$ fast-growing sequences unlimited number of times (for variable input length), this corresponds to an infinite sequence of level $k + 1$. However, to make the (arguably) finitistic definition as concrete as possible (and to allow a more fine-grained hierarchy of notions), we limit ourselves to a single application (of the producer of the sequences), and to correctness with respect to an ordinary Turing machine using that application.

**Restrictions on Length**
(Disclaimer: Statements here were not formally verified.)
   In addition to considering sufficiently long sequences, we can also consider fast-growing sequences $g$ with specific length conditions. Given a recursive well-founded relation '≺', a natural corresponding length condition is the following: There is no $x_0 \succ x_1 \succ \ldots \succ x_n$ with $x_i < g(i)$ and $n + 1$ being the length of $g$. Computability for sufficiently fast-growing $g$ with this length condition equals $0^{(\alpha+1)}$ (or $0^{(\alpha)}$ for finite $\alpha$) where $\alpha$ is the rank of '$\prec$'.

For fixed finite $k$, $g$ having length $k$ corresponds to ordinal $k$ for '$\prec$'. If $h$ is a fixed recursive function with $\lim_{n\to\infty} h(n) = \infty$, then $g$ having length $h(g(0))$ corresponds to (computability using) $0^{(\omega+1)}$. The independence of expressiveness on $h$ may appear counterintuitive — if we increase the length of $g$ by 1, we can compute the halting problem for appropriate machines that use the shorter $g$ — but we cannot determine which machines are valid (in the sense of being independent of the choice of sequence), hence halting reversal does not increase expressiveness. Going further, length $h(g(k))$ corresponds to $0^{(\omega+k+1)}$, length $g(h(g(k))) - 0^{(\omega 2+k+1)}$, $g^m(h(g(k))) - 0^{(\omega(m+1)+k+1)}$ (superscript indicates repetition), $g^{g^m(h(g(k)))}(0) - 0^{(\omega^2+\omega m+k+1)}$, and similarly using nested superscripts and ordinals below $\omega^3$.

**Axiomatization**

It is insightful to axiomatize key properties of level $k$ sequences. A recurrent theme is that axioms using infinity give us symmetries and that using 'sufficiently large' (and constructs that build on it), we can axiomatize finite versions of these symmetries.

**Using nonstandard numbers:** Intuitively, a sufficiently long level $k$ sequence (or another construct using Definition 4.2) can often be treated as if it were a definite partially unknown object, and we can formalize this analogously to Robinson nonstandard analysis. Add a predicate symbol $I$ for small numbers (which will include all standard numbers), an axiom that $I$ is a proper cut, and a schema (over definitions) that if $x$ is small and $y$ is definable from $x$ then $y$ is small. In Definition 4.2, augment the directed system so that it is a filter. Given a language L, a sufficiently good example will be an example that satisfies all notions with a 'small' definition in L. An alternative is to work inside L augmented with $I$, add a constant symbol for some sufficiently good example, and add appropriate axioms that are true for every sufficiently good example.

If one wants an axiomatization without nonstandard numbers, infinite structures, or a function symbol $f$:$n$→(an example sufficiently good relative to $n$), one can do as follows. Add a symbol for a sufficiently good example and (optionally) a predicate for small numbers. Identify a class $F$ of functions that we know are independent of the example provided that the example is sufficiently good relative to the argument of the function. Add axioms: (1) (schema, $f$ in $F$) $f(0)$ is small, (2) (schema, $f$ in $F$) for small $x$, $f(x)$ is well-behaved, and perhaps other schemas of that nature. Prove theorems as usual, but when interpreting the results, assume that the schema works only up to some sufficiently large code "f", and is thus consistent with a particular example.

Variants of the finitary treatment above include having two notions of "small", so that invariant small-1-length definitions only define numbers that are small-2, or one can have more notions of "small", and partially simulate "small" being a cut. Alternatively, one can avoid "small" as a predicate, and instead interpret "small" using the first element (or equivalent) of the sequence (and the second element for the second notion of "small", and so on), and remove these elements from the sequence (or do an analogous operation).

We now attempt to axiomatize nonlocal level $k$ sequences. We will use arithmetic augmented with a symbol meant to represent a level $k$ sequence, and the intent is to capture key properties of the sequence without using infinite sets. The axiomatization will be generic, allowing both finite and infinite sequences, with an added condition to choose the type and length. Its strength and consistency are unclear. For the finite version, the axioms are $\Delta_2^0$ above the sequence, and they might have a $\Pi_1^0$ strengthening by requiring that a witness for $\varphi$ does not increase by the reduction in $q$. If a bounded axiomatization is desired (so that satisfiability of a finite fragment is $\Sigma_1^0$), use φ for which (syntactically) any witness must be bounded above $A$, and (if desired) use a schema over cut-off points to give the axiomatization a non-arbitrary strength.

The sequences will be coded as predicates and strictly monotonic. A *bounded shift* of $A$ above $p$ is a strictly monotonic $A'$ with $A(x+p) \leq A'(x) \leq A(A(x+p))$. $A'(x)$ exists iff $A(A(x+p))$ exists. The idea is that a bounded shift will preserve the desired structure of $A$ while erasing (or be capable of erasing) detectable particularities of the exact choice of numbers.

Let $\varphi$ be a reasonable universal $\Sigma_1^0$-formula. We say that $\varphi(p,q)$ is *level $i$ invariantly true* iff $\varphi(p,q)$ holds using $A_1, \ldots, A_{i-1}, A'_i, \ldots, A'_k$ for every combination of bounded shifts above $\max(p,q)$. We formalize this as a $\Sigma_1^0$-formula by requiring a bound on the witness that works for all shifts (by König's lemma, this is correct).

**Axioms:**
* Base system such as EFA or PRA or $\Sigma_1^0$-PA or PA.
* Basic properties of $A$:
   - $A_i \cap n$ exists as a number
   - $A_i$ is strictly monotonic, and if $A_i(n)$ exists and $m < n$, then $A_i(m)$ exists.
   - If $i < j$, then $\max(A_j) \leq \max(A_i)$ (and if $A_j$ is infinite, then so is $A_i$).
* Formalization of sufficiently large:
   - If it exists, the least $q$ such that $\varphi(p,q)$ is level $i$ invariantly true is below $A_i(p)$ (if $A_i(p)$ exists).
   - (schema over $p$) (optional, or can be as long as desired) If it exists, the least $q$ such that $\varphi(p,q)$ is level $i$ invariantly true is below $A_i(0)$.

Here is one choice for the finite sufficiently long condition. Different choices give different strengths. For level 1, the strength here appears to correspond to $\Pi_1^1 - \mathrm{CA}_0$.
* (optional) Each $A_i$ is finite.
* (schema, or a particular choice of $n$) There is a bounded shift (above 0) $(A'_1, \ldots, A'_k)$ such that $(A'_1, \ldots, A'_k, A'_{k+1})$ satisfy the axioms for level $k+1$ sequences for some $A'_{k+1}$ such that $A'^n_{k+1}(0)$ exists ($n$ indicates iteration).

*Note:* For infinite $A$, simply asserting their infinity should give a good theory. The second statement can still be used for bounded shifts of sufficiently long initial segments of $A$.

If true and not derivable from more basic axioms, we also want to add the following: If a $\Sigma_1^0$ statement $p$ is true for some bounded shift $(A'_1, \ldots, A'_{k-1})$ above $p$, then for every $p' \geq p$, $p$ is true for some bounded shift $(A'_2, \ldots, A'_k)$ above $p'$ ($p$ codes the statement; the statement depends on $k - 1$ second order parameters). This applies to the infinite version; and if $A_{k+2}(p')$ exists, also to the finite version. The added condition ensures that $A_{k+1}$ is sufficiently long. Note that a $\Sigma^0_1$ statement cannot directly test whether an end of $A$ is reached. Also, without the statement, it is not clear whether different levels of $A$ agree on well-foundness. The statement can be supplemented with a schema over $p$ using bounded shifts above 0 and $p'$ (and place of $p$ and $p'$). For an infinite level 1 sequence $A$, we may also want to state (if not already provable) that well-foundness of a $\Sigma^0_1$ relation that is invariant above $p$ is unchanged by using a definable (or restricted definable) infinite subsequence of $A$ in place of $A$. This also applies to the highest level of a higher level sequence if trying to get the most out if it.

# 5 $\Sigma^1_2$ Truth

We now give an arguably finitistic definition of $\Sigma^1_2$ truth. Essentially, this section informally describes in detail a variation on level 2 estimators of the next section.

In the computational model, along with input, we are given a finite set of finite sequences. Each sequence is sufficiently long relative to itself, or relative to its rate of growth. (Optional) If $A$ is sufficiently long and $|A| \leq |B|$ and $\forall n < |A|\, B(n) \leq A(n)$, then $B$ may be considered sufficiently long. (Optional) For the purposes of determining sufficient completeness (below), arbitrary conditions may be imposed on the sequences as long as they do not impose an upper bound on the rate of growth (as used here, an upper bound would be a function $f$ such that for every permitted sequence $A$, $\exists n < |A|\, A(n) < f(n)$). Given these conditions, the set is sufficiently complete for the input size. Completeness is preserved under addition of new sequences or decrements to input size. We are not given the conditions on the sequences, but regardless of what they are, the algorithm below is guaranteed to be correct.

The algorithm aims to check whether the sentence encoded in the input has a transitive model. It enumerates potential initial segments (whose length equals the length of the longest sequence in the set of sequences we were given) of the truth predicate for the language with a distinguished symbol for a well-ordering of the elements. Each segment is checked for a quick (relative to its length) inconsistency with the sentence. If not found inconsistent, it is converted into a finite quasi-model, with elements coded by sentences of the right syntactic form: $\exists x\, P(x)$ codes the least such $x$ if true, and 0 if false. Then, for every sequence $A$ (in the set of sequences we were given) with $A(0)$ larger than the input size, it tries to form a descending path $(a, b, c, \ldots)$ of elements with $a < A(1), b < A(2), c < A(3), \ldots$. If the path reaches length $|A| - 1$, the truth predicate is rejected; otherwise it is accepted by the sequence. The algorithm returns 'yes' iff a truth predicate segment is accepted by all sequences.

To see that the computation is correct, note that if a sentence has a transitive model, then it has a transitive model of bounded complexity, and hence without sufficiently long descending sequences of the required type. If it does not have a transitive model, then in the complete binary tree of candidate truth predicate initial segments, every path will have an infinite descending sequence (or inconsistency). Thus, regardless of how strictly "sufficiently long" is defined, every path will have sufficiently long descending sequences; furthermore, excluding all sequences that are above an infinite path would be an upper bound on the rate of growth (as used above), which is prohibited. By König's lemma, a finite set of such sequences can cover all paths, and thus a sufficiently complete finite set will suffice.

The construction has an infinite version, namely a predicate or a tree on finite sequences of natural numbers. Consider a well-founded tree of such sequences (ordered by extension) such that every non-leaf node can be extended by any natural number, and if a sequence is a leaf node, then it is sufficiently long relative to itself (or relative to its rate of growth; the difference does not affect the expressive power) — or any tree that is a well-founded extension of such a tree. Using such a tree, recognizability is $\Pi^1_2$ complete. (However, $\Sigma^1_2$ may be a better analogue of being recursively enumerable: The relevant parameter is not the halting time but the ordinal depth of the visited portion of the tree for nonhalting runs, which corresponds to the complexity of the least witness for a $\Sigma^1_2$ statement.) In one direction, the above computational procedure works, except that if the theory has a well-founded model, we do not halt. In the other direction, the machine halts iff for every well-founded tree $T$ for every non-halting computation corresponding to an extension $S$ of $T$, $S$ has an infinite path, which is $\Pi^1_2$.

We can extend definability by using an infinite sequence of numbers growing sufficiently fast relative to the tree. Recognizability should then correspond to having a $\Pi^1_2$ monotonic inductive definition. A stronger extension (likely still $\Delta^1_3$) is to use an $n$-tuple of trees, with sequences in each tree sufficiently long relative to themselves, and to every previous tree.

# 6 Second Order Arithmetic

**Vague Definition 6.1:**
A *level 0 estimator* is a pair of natural numbers $(a, b)$ with $a$ sufficiently large, and $b$ sufficiently large relative to $a$.
A *level n+1 estimator* is a sufficiently complete finite set of level $n$ estimators, where 'sufficient completeness' is defined such that it does not prevent the notion of level $n$ estimators from being arbitrarily strict.
(Sufficiency depends on the problem/calculation considered.)

*Notes:*
* Given a finite set of problems, a single notion of sufficiency works. Given a countable list of problems, a single notion of sufficiency works, except that minimum $a$ depends on the problem size.
* Without affecting the expressive power, instead of 0-estimators, we could have started with level -1 estimators being sufficiently large numbers (or even every number being a sufficiently good level -2 estimator).
* A variation that leads to the same expressiveness for $n > 0$ estimators is to drop the condition on $a$; $b$ will be sufficiently large relative to $a$ and the problem considered.
* An example of the "arbitrarily strict" qualifier is that even if we admit $(a, b)$ as a level 0 estimator, there are sufficiently complete level 1 estimators with every included $(c, d)$ satisfying $d > 2b$.

The above definition is vague but the following is well-defined regardless of how the vagueness is resolved. (As written, Definition 6.2 is well-defined relative to the choice of an estimator and (for bad estimators) the choice of 'standard rules', but given Theorem 6.6, neither is relevant if we use Definitions 6.3 and 6.4.)

**Definition 6.2:** Computation of second order arithmetical truth:
* Use standard rules to convert the sentence (or its negation) to $QX_1QX_2\ldots \exists X_n \forall x \exists y > x\, P$ where $n > 0$, each $Q$ is a quantifier, and $P$ is a bounded quantifier formula that does not use $x$ (and $P$ can even be chosen to be polynomial time computable), and each $X_i$ is a set of natural numbers.
* Pick any (it does not matter which one) level $n$ estimator on problem sizes at least as large as the sentence. The sentence is true iff it passes the estimator.
* A formula $\forall x \exists y > x\, P(y, \ldots)$ ($P$ as above) passes a level zero estimator $(a, b)$ iff $\exists y (a \leq y < b) P(y, \ldots)$.
* A formula $\forall X \varphi(X)$ passes an estimator iff for every binary sequence $X$, $\varphi(X)$ passes at least one of the estimators included.
* A formula $\exists X \varphi(X)$ passes an estimator iff for some binary sequence $X$, $\varphi(X)$ passes all of the estimators included.
(Note: Because $P$ is bounded quantifier and an estimator is finite, the search over $X$ is also finite).

To intuitively see that the description is correct, if a $\Sigma^1_n$ statement is true, then it has a witness (for the outermost existential quantifier) of bounded complexity. Level $n-1$ estimators will work correctly with that witness and hence the sentence. On the other hand, if a $\Sigma^1_n$ statement is false, then every potential witness will be intercepted by some level $n-1$ estimator, and hence the falsehood will be witnessed by a sufficiently complete set of such estimators. By induction on $n$, the method is correct.

To prove correctness, we first have to formalize what we want to prove. If we treat 'sufficiently' as a parameter in the definition, we naturally get the following:
**Definition 6.3:**
* A notion $R$ of 0-estimators is valid iff there are $a_{\min} \in \mathbb{N}$ and $f \in \mathbb{N}^{\mathbb{N}}$ such that $R(a, b) \Leftrightarrow a \geq a_{\min} \wedge b \geq f(a)$.
* A notion $R$ of $n+1$ estimators is valid iff
    (1) Every element of $R$ is finite. (A reasonable variation (immaterial for below): Also require elements to be non-empty.)
    (2) $\cup R$ is valid and $x \in R \wedge y \in \cup R \Rightarrow x \cup \{y\} \in R$. (Equivalently, there is a valid notion $T$ of $n$-estimators such that every element of $R$ consists of elements of $T$, and adding an element of $T$ to an element of $R$ results in an element of $R$.)
    (3) For every valid notion $S$ of $n$-estimators, there is an element of $R$ consisting of elements of $S$, i.e. $\exists r \in R\, r \subset S$.

**Definition 6.4:** A sentence (in second order arithmetic or another well-defined language) using a (sufficiently good) $n$-estimator as a parameter is well-defined iff there is a valid notion $R$ such that all elements of $R$ agree on whether the sentence is true, which is then the truth value of the sentence.
*Note:* This is simply Definition 4.2(c) applied to valid notions of estimators.

The truth value of a well-defined sentence is unambiguous:

**Proposition 6.5:**
**(a)** Every valid notion of $n$-estimators is non-empty.

**(b)** Given two valid notions $R$ and $S$ of $n$-estimators, $R \cap S$ is also a valid notion.
**Proof:**
**(a)** This follows from 6.3(3) and existence of the trivial notion of $n$-estimators that includes everything of the right syntactic form.
**(b)** Assuming $n > 0$ ($n = 0$ is easy), note that $\cup R$ and $\cup S$ are valid. Assume that the proposition applies for $n - 1$, and let $T$ be an arbitrary notion of $n - 1$ estimators. Thus, $T' = T \cap \cup R \cap \cup S$ is valid. There is $x \in R$ and $y \in S$ such that $x \subset T'$ and $y \subset T'$, and for every such $x, y$, we have $x \cup y \in R \cap S$, which satisfies (3) in the definition. Using a similar argument, $\cup R \cap \cup S = \cup (R \cap S)$, and (2) follows.

**Theorem 6.6:** The above notion/computation of second order arithmetical truth is well-defined and agrees with the actual truth.
**Proof:** We will use a conservative extension with nonstandard numbers (see discussion in Section 4). Augment the language with a well-ordering of real numbers (this simplifies the proof if the base theory does not prove projective uniformization). Let $I$ be a proper cut of natural numbers closed under exponentiation. An $X$-estimator will be an estimator that is valid for every valid notion definable in second order arithmetic (with the well-ordering) using $X$ as a parameter, with the definition in $I$. An estimator is an $X$-estimator for the empty $X$. We show by induction on $k$ that the computation works for $k$-$X$ estimators for formulas in $I$ that use $X$ as a parameter. $k = 0$ is trivial.
If $\exists Y \varphi(X, Y)$ is true, then $\varphi(X, Y)$ holds for some $Y$ of bounded complexity relative to $X$, and hence every $k$-$X$-estimator will handle $\varphi(X, Y)$ correctly.
If $\exists Y \varphi(X, Y)$ is false, then for every $Y$, every $k$-$X$-$Y$ estimator will reject $\varphi(X, Y)$. Thus, for every valid notion of $k$-estimators, using König's lemma, there is a finite set $S$ of $k$-estimators such that for every $Y$ there is $s \in S$ that rejects $\varphi(X, Y)$. Therefore, every $k+1$-$X$-estimator rejects $\exists Y \varphi(X, Y)$.
$\forall Y \varphi(X, Y)$ is analogous.
*Note:* Using $I$ spared us some complexity of tracking permitted problem sizes.

If sufficiently complete did not obey the "arbitrarily strict" (for previous levels) qualification and we used a definition like Definition 3.1, the expressiveness using $k$-estimators ($k > 0$) should collapse to that of being hyperarithmetic in the $k$th hyperjump of 0. (Essentially, enough level 0 estimators give us a sufficiently fast-growing sequence $S$, and levels 1-k estimators would give us $k$ special points (each sufficiently large relative to the previous points and the growth rate of $S$), including the final point of $S$.)

For readability, we sometimes use $X$-estimator to mean an estimator sufficiently good relative to $X$. The reason is that we often intuitively (and sometimes formally) reason as if 'sufficiently good' has a concrete meaning, but that only works if 'sufficiently good' is parameterized by the free variables of the formula. An example (giving a key property of estimators) is that for every continuous predicate $P$ on $(2^{\mathbb{N}}, \mathbb{N})$ and $P$-estimators $s$ and $t$, $\forall X \exists s' \subset s \, P(X, s') \Leftrightarrow \forall X \exists t' \subset t \, P(X, t')$.

The estimators have an infinite version: A predicate that codes a sufficiently strict valid notion of $n$-estimators. Recognizability for machines using such predicates (where the machine must work for every sufficiently strict notion of $n$-estimators) appears to equal $\Pi^1_{n+1}$. Also, even without the condition on $a$ for 0-estimators, a single predicate suffices for all input sizes: Given input size $s$, the machine can pick $n$

-estimators such that for every included level 0 estimator $(a, b)$, $b > s$, and then ignore 0-estimators $(a', b')$ with $a' < s$.

This completes our description at the level of second order arithmetic. What remains is finitistically studying which axioms to include in the system, and relating infinite statements to natural statements about finite structures. If, for example, projective determinacy for the results of the estimations will be deemed a finitistically acceptable axiom schema, we will have a finitistic consistency proof of ZFC.

We take a key step in this direction by interpreting projective determinacy in terms of finite games and a basic interaction of adjacent levels of estimators. A possible intuitive interpretation is that completeness for level $k + 2$ estimators does not negate completeness for level $k + 1$ estimators.

**Theorem 6.7:**
**(a)** Lightface projective determinacy holds iff for every $k > 0$ and every primitive recursive predicate $P(n, X, Y)$ ($X,Y$ and (below) $Z$ are predicates on natural numbers), for some (equivalently, every) level $k + 2$ estimator $s$ (sufficiently good relative to $P$), one of the players has a winning strategy in the following game: player I - $X_0$, player II - $Y_0$, player I $X_1$, player II $Y_1$, ...
player I wins iff $\forall s_{k+1} \in s \exists s_k \in s_{k+1} P(s_k, X, Y)$
player II wins iff $\forall s_{k+1} \in s \exists s_k \in s_{k+1} \neg P(s_k, X, Y)$
**(b)** Projective determinacy holds iff for every $k > 0$ and every $P(n, X, Y, Z)$ (as above), for some (equivalently, every) level $k + 3$ estimator $t$ (sufficiently good relative to $P$), for every $Z$, there is $s \in t$ such that one of the players has a winning strategy in the above game (using $P(s_k, X, Y, Z)$).
*Notes:*
* By virtue of the bounds in $P$ (which can be extracted from a primitive recursive definition of $P$), $X$, $Y$, and $Z$ can be treated as bounded sets, and the game is finite.
* For every $n$, there is a single well-behaved $P$ sufficient for $\Sigma_n^1$ determinacy (both lightface and boldface): Player I chooses a $\Pi_n^1$ game and wins if he wins that game and player II does not play a code for a winning strategy for player I. Also, given an estimator $s$ (or $t$) testing who wins (and whether the determinacy holds) is primitive recursive.
* Instead of "primitive recursive", we could have used "elementary time", or in the other direction, continuous (and for (a), continuous and definable in second order arithmetic). König's lemma is used implicitly in the proof.
* We do not require for $P$ to be well-behaved. $P$ is well-behaved iff for every $(X, Y)$, $P(n, X, Y)$ is independent of $n$ if $n$ codes a sufficiently good (relative to $P, X, Y$) $k$-estimator. One can show that an equivalent formalization of well-behaved is that for some (equivalently, every) estimator $(k + 2)$-$P$-estimator $s$,
$\forall X, Y \exists s_{k+1} \in s \forall a \in s_k \forall b \in s_k\ P(a, X, Y) = P(b, X, Y)$.
**Proof:**
**(a)** We show that the first player can win the finite game iff he can win the infinite game corresponding to $P$, with the ambiguity in the construction of $P$ resolved in favor of the first player: $P_1(X, Y)$ iff for every notion of $k$-estimators there is a $k$-estimator $n$ with $P(n, X, Y)$, equivalently iff for every (equivalently, some) $(k + 1)$-$P$-$X$-$Y$-estimator $s_{k+1}$, $\exists s_k \in s_{k+1} P(s_k, X, Y)$. Combined with the same result for the second player, and given the complexity of the formulas, this leads to the equivalence with the determinacy.

Let $G$ be the function converting (player I strategy, player II strategy) into the corresponding plays. Note that $G$ is the same for both the finite and infinite game. The first player has a winning strategy in the infinite game iff
$\exists S \, \forall T \, \exists (k\text{-}S\text{-}T \text{ estimator } s_k) \, P(s_k, G(S, T))$.
Fix an $S$ such that $S$ is winning if player I can win. Thus, we can assume that all estimators are sufficiently good relative to $S$, and the winning condition becomes $\forall T \, \exists (k\text{-}T \text{ estimator } s_k) \, P(s_k, G(S, T))$. This implies that there is a finite set $s_{k+1}$ of level $k$ estimators such that $\forall T \, \exists s_k \in s_{k+1} P(s_k, G(S, T))$, and therefore every $k+1$-estimator $s_{k+1}$ will work, and the first player wins the finite game by using $S$. Conversely, if the first player has a winning strategy in the finite game, fix a notion $k+1$-$P$-estimators. Build an infinite tree of first player strategies (ordered by extension) for the first $n$ moves that do not lose against any of the estimators (here, an incomplete game is not a loss), and pick an infinite path (and strategy) $S$. The first player wins the infinite game using $S$ because for every $T$ he wins the finite game corresponding to a $(k+2)$-$S$-$T$ estimator.
**(b)** The proof is analogous to *(a)* with an extra quantifier pair to handle $Z$, or alternatively, easily follows from the relativized version of *(a)*.

**Axiomatization**

While the above theorem is an equivalence and not just an equiconsistency, we also want an equiconsistency with a simple intuitive axiomatization of estimators that does not use infinite sets. We expect that the following axiomatization works.

A $k$-$n$-estimator will be an $k$-estimator sufficiently good relative to $n$. We idealize 'sufficiently good' with a schema such that if the schema is cut off at any finite point, there are particular notions consistent with the axioms. An alternative would be to explicitly adjust $n$ based on the formula size (and in the axiom 7 and definition of explicitly invariant below, use a number coding the subformula, the free variables, and $k$ in place of $\max(x_2, \ldots, x_i)$).

We start with arithmetic, and add a predicate for $k$-$n$-estimators, with one predicate symbol per $k$. An alternative to a predicate would be a function: $n \to$ example of $k$-$n$-estimator, and application of basic closure operations would lead to a predicate consistent with the axioms. Also, with some added complexity, one can axiomatize a single example of a $k$-$n$-estimator and using minimal arithmetical axioms.

*Definition:* A formula is explicitly invariant if all of its uses of $k$-$n$-estimators are in the form $\exists s$ ($s$ is $k$-$n$-estimator) ... where $n$ is the maximum of the free variables (not including $s$) at that location of the formula, and $k$ is a numeral. Finite sets will be coded by numbers (using any reasonable coding).
*Note:* The intent is that by virtue of its syntactic form, the formula is independent of how strictly we construe 'sufficiently complete'.

**Language:** arithmetic (using natural numbers), with predicates for $k$-$n$-estimators, with one predicate symbol per $k \geq -2$.
**Axioms:**
1. Basic arithmetical axioms.
2. Induction for all formulas. [note: limited induction likely also works]
3. Every number is a -2-$n$-estimator.
4. A $k+1$-$n$ estimator is a finite set of $k$-$n$-estimators.
5. Adding a $k$-$n$-estimator as an element to a $k+1$-$n$-estimator results in a $k+1$-$n$

-estimator.
6. Every $k$-$n+1$-estimator is a $k$-$n$-estimator.
7. (schema over explicitly invariant $\varphi(x_1, \ldots, x_i)$) whether $\exists s' \subset s\, \varphi(s', x_2, \ldots, x_i)$ holds is independent of the choice of the $k$-$\max(x_2, \ldots, x_i)$-estimator $s$. [note: symbol $\subset$ permits equality, but that is immaterial here]
8. For every $n$, a $k$-$n$-estimator exists.

Axiom 7 (combined with axioms 4 and 5) simply asserts that $s$ is sufficiently complete as tested using $\varphi$ and $x_2, \ldots, x_i$.

The axioms are consistent and sound in that every finite fragment is satisfied using every sufficiently strict function: $n \to$ valid notion of $k$-$n$-estimators (valid as in Definition 6.3; $k$ is the maximum $k$ used in the fragment).

*Variations on Explicitly Invariant:* The class of explicitly invariant formulas can be broadened somewhat: With the axiomatization remaining sound, max can be replaced (including in axiom 7) by $f(\max(\ldots))$ (using a schema over $f$) where $f$ is explicitly invariant and is defined without using free variables and $\lim_{n\to\infty} f(n) = \infty$. If the notion of $k$-$(n+1)$-estimators is sufficiently strict relative to $k$-$n$ estimators, the class can be broadened further still. It is unclear if we need any restriction for $\varphi(x_1, \ldots, x_i)$ other than (1) every use of estimators is $\exists s$ ($s$ is $k'$-$n$-estimator) where $n \geq f(\max(x_1, \ldots, x_i))$ ($f$ as above; also note that the variables that are free in the subformula of φ but not free in φ are unrestricted) and (2) a positive $\exists s$ ($s$ is $k$-$n$-estimator) … cannot appear inside a negative $\exists s$ ($s$ is $k'$-$n$-estimator) … and vice versa.

We conjecture that the axioms are equiconsistent with second order arithmetic $Z_2$, and that the addition corresponding to the above theorem (using a schema over $P$, or a quantifier over codes and $k$-code(P)-estimators, or even a single universal instance) is equiconsistent with $Z_2$ + lightface projective determinacy for *(a)* and $Z_2$ + projective determinacy for *(b)*.

# 7 Beyond Second Order Arithmetic

## 7.1 Estimators with Transfinite Levels

To go further, we can define transfinite levels of estimators. For an ordinal $\alpha > 0$ (alternatively, for consistency with $n$-estimators, limit $\alpha$), an $\alpha$-estimator is a sufficiently complete collection of estimators of lower levels, with α being marked in the estimator, and where 'sufficiently' complete is defined so as not to preclude notions of lower level estimators being arbitrarily strict. Level 0 estimators are as defined above (Definition 6.3 or an alternative).

**Definition 7.1.1:**
A notion $R$ of $\alpha$-estimators ($\alpha > 0$) is valid iff:
(1) Every element of $R$ is finite.
(2) $\cup R$ is a valid notion of $<\alpha$-estimators.
(3) $x \in R \wedge y \in \cup R \Rightarrow x \cup \{y\} \in R$
(4) For every valid notion $S$ of $<\alpha$-estimators, there is an element of $R$ all of whose elements are elements of $S$.

A notion $R$ of $< \alpha$-estimators is valid iff:
(5) For every $\beta < \alpha$, the set of $\beta$-estimators in $R$, $R_\beta$, forms a valid notion.
(6) (Coherence) For every $\beta < \alpha$, $R|\beta = \cup R_\beta$. (Here $\beta > 0$ and $R|\beta = \cup_{\gamma<\beta} R_\gamma$. Informally, different levels of $R$ agree on which estimators of lower levels are valid.)

**Proposition 7.1.2:**
**(a)** For every ordinal $\alpha$, there is a valid notion of $\alpha$-estimators, and every valid notion is non-empty.
**(b)** For valid notions $R$ and $S$ of $\alpha$-estimators, $R \cap S$ is valid. Similarly, intersection of two valid notions of $< \alpha$-estimators is valid. Also, for $\alpha > 0$, $\cup(R \cap S) = \cup R \cap \cup S$.
**Proof:**
**(a)** This follows from the existence of the trivial notion that includes every set of the right syntactic form. (Also, any valid notion of $< \alpha$-estimators has a trivial corresponding notion of $\alpha$-estimators.)
**(b)** Assume this holds for smaller $\alpha$ (and $\alpha = 0$ is immediate). Let $R'$ be $\cup R$ and $S'$ be $\cup S$. Verify that $T' = R' \cap S'$ is a valid notion of $< \alpha$-estimators: (5) $T'_\beta = R'_\beta \cap S'_\beta$, which is valid, and (6) $T'|\beta = R'|\beta \cap S'|\beta = \cup R'_\beta \cap \cup S'_\beta = \cup(R'_\beta \cap S'_\beta) = \cup T'_\beta$. Thus, there are $r \subset T'$ with $r \in R$ and $s \subset T'$ with $s \in S$, so $r \cup s \in R \cap S$, and (using (3) for $R$ and $S$) every element of T' is in $\cup(R \cap S)$, so $\cup(R \cap S) = \cup R \cap \cup S$. Finally, (4) holds because for a valid notion $U'$ of $< \alpha$-estimators, $T' \cap U'$ is also valid.

Thus, using Definition 4.2, we can meaningfully speak of 'for every sufficiently good $\alpha$-estimator'.

Definability of reals using estimators of levels below $\omega(\alpha + 1)$ appears to equal to definability in $L_\alpha(\mathbb{R})$, at least when $\alpha$ is not too large. Note that for a definition using $\alpha$-estimators not to (overtly) use infinity, we would need a coding of ordinals up to $\alpha$ by integers.

**Estimators with Rational Levels**

Countable ordinals can be embedded into rationals, and one extension of transfinite level estimators is to use estimators with rational levels. Level 0 estimators will be as before, and a higher level estimator will consist of the level and a sufficiently complete finite set of lower level estimators, but the treatment of 'sufficiently complete' is more subtle. Rational numbers are not well-ordered, but we work around that by in effect defining estimators for sufficiently complete (unspecified) well-ordered subsets of nonnegative rationals.

Given a notion $R$ of $\alpha$-estimators and a zero-preserving order-preserving injection $f : \alpha \to \beta$, there is a natural translation of $R$ to $\beta$, $R_f$:
$\gamma$ not in the range of $f$ - every finite subset of $R_f|\gamma$ is included (as a level $\gamma$ estimator in $R_f$).
$\gamma$ in the range of $f$ - a finite subset $r$ of $R_f|\gamma$ is included if after recursively deleting from $r$ all estimators of levels that are not in the range of $f$, removing duplicates, and applying the inverse of $f$ on the levels, we get an element of $R$.

The above definition works even if $\beta$ is a non-well-founded linear order, and in particular for $f : \alpha \to \mathbb{Q}$. The set of notions of estimators derived from a notion of $\alpha$-estimator and $f : \alpha \to \mathbb{Q}$ (variable countable $\alpha$ and $f$) forms a directed system under strictness, which allows Definition 4.2 to work.

The expressiveness for computations that use sufficiently good estimators with rational levels probably equals to $\Delta_1(\{\mathbb{R}\})$ definability in $L_\alpha(\mathbb{R})$ for the least $\alpha$ such that $\alpha$ is not $\Delta_1(\{\mathbb{R}\})$ definable in $L_\alpha(\mathbb{R})$ (and thus $\alpha$ is large but countable here).

If the system used finite ordinals instead of countable ordinals, the expressiveness appears to equal the Turing degree of second order arithmetic. The difference is that (with finite ordinals) if we use (for example) a level 1 estimator to get estimators of level 1/2, 2/3, 3/4, 4/5, ..., we still get enough levels of estimators, but there is no assurance that an estimator of level $(n-1)/n$ is sufficiently complete relative to $n$ for large $n$ (or even non-empty, and we do not know when to stop considering higher level estimators).

A caveat (even for the full system above) is that if we use (for example) a level 1 estimator to get estimators of level 2/3, 3/5, 4/7, ..., then for a sufficiently large $n$, we may get an empty (or otherwise small) estimator for level $n/(2n-1)$. Of course, if $n$ is not too large (or if the completeness is relative to $n$), all included level $n/(2n-1)$ estimators would still be (essentially) sufficiently complete. There might be interesting notions of estimators with rational levels that handle such descending sequences differently (or other interesting notions of estimators with rational levels in general). Perhaps, such notions can be used to compute game quantifiers (and strategies) for games of countable length.

## 7.2 Using Pattern Enumerators

A (not necessarily finitistic) notion related to that of a sufficiently large Turing degree (or more closely, fast-growing infinite sequence of Turing degrees) is that of a sufficiently closed pattern enumerator.

**Vague Definition 7.2.1:** A sufficiently closed pattern enumerator is a function $f : \mathbb{N} \to 2^\mathbb{N}$ such that
- $f(i)(n)$ is character $n$ in the infinite binary string corresponding to pattern $i$. Numbering of patterns is permitted to be arbitrary.
- if $x \in 2^\mathbb{N}$ follows a sufficiently simple pattern relative to $f(i)$, then $x$ is in the range of $f$.
- (optional for our purposes) if $x$ and $y$ are in the range of $f$, then so is their join $x_0, y_0, x_1, y_1, \ldots$.

For example, for some $i$, $f(i)(n) = 0$ for all $n$, but the least such $i$ may be arbitrarily large. Also, without the optional condition, enumerators $A$ and $B$ can be trivially combined into an enumerator of patterns included in $A$ or $B$.

As before, we can use Definition 4.2 for unambiguous usage of sufficiently closed pattern enumerators. A valid notion of 'sufficiently simple' is simply a function that for every element of $2^\mathbb{N}$ assigns a countable subset $2^\mathbb{N}$, with stricter notions assigning more inclusive subsets, and pattern enumerators are defined accordingly.

Let $A_1, A_2, \ldots, A_n$ be sufficiently closed pattern enumerators such that $A_{i+1}$ includes $A_i$ as a pattern, and $S$ be a sufficiently fast-growing sequence relative to the enumerators.

Assuming sufficient determinacy, we can recognize the set of winning positions of games of integers of length $\omega^2 n$ with $\Sigma_1^0$ payoff as follows. Replace the first $\omega$ moves of player I by a strategy $p_1$ for these moves where $p_1$ (chosen by player I) is in $A_1$, replace the first $\omega$ moves of player II by strategy $s_1$, such that $p_1$ join $s_1$ is in $A_1$, and so on replacing the first $\omega^2$ moves with $p_1, s_1, p_2, s_2, \ldots$ and reducing game length to $\omega^2(n-1)$ (if $n > 1$). We might not be able to go further with $A_1$ since $A_1$ does not recognize its behavior as a pattern, but we can use $A_2$ for the first $\omega^2$ moves of the modified game and so on, converting the game to one of length $\omega$. We then constrain the move values using $S$, and if there is a strategy (for the first player) that wins at a finite stage, halt with 'accept'. We conjecture the converse to hold as well (assuming sufficient determinacy).

*Notes:*
* Without a sufficiently fast-growing sequence $S$, an unspecified sufficiently closed pattern enumerator is essentially useless to a Turing machine since it may be arbitrary up to an arbitrarily large finite point.
* Sufficiently closed pattern enumerators have a number of equivalent (in terms of expressiveness) definitions. For example, we may view them as pattern detectors with $f'(s)$ being the least $n$ such that $f(n)$ starts with $s$.
* To go further, we can use a transfinite sequence of pattern enumerators, with each enumerator including the sequence of enumerators of lower levels (combined into a binary sequence) as a pattern.
* If we used 'sufficiently simple pattern' instead 'sufficiently simple relative pattern relative to $f(n)$', we should get the same expressiveness by using $\omega\alpha$ enumerators in place of $\alpha$, while requiring each enumerator to include all sufficiently simple patterns with respect to the join of the lower level enumerators.
* Like fast-growing sequences, the enumerators (and chains of them with the fast-growing sequences), have a finite version; simply cut off the chain (coded as an infinite string) at a sufficiently large number relative to the chain and the problem size.

*Weak versions and computational complexity:* Let an enumerator be constant but unknown. Recognizability for polynomial time computation equals NP (no time limit on nonaccepting computations; assumes binary access to indices; unary access corresponds to unary NP under many-one polynomial time reducibility). If we are also given a point such that every sufficiently simple pattern comes before it, P-time decidability equals PSPACE, and if we also have access to randomness, P-time (probabilistic) decidability captures $E^{NP}$ (and hence NEXP), but is within $\Delta_2$-EXP. If we are only required to be correct for sufficiently large instances (without using the large point above), decidability appears to equal $\Delta_1^1$ (with recognizability $\Pi_1^1$). Arithmetic definability using the enumerator (equivalently, computability using a fixed length sufficiently fast-growing sequence) equals definability in second order arithmetic. However, if the enumerator is not required to be closed relative to every pattern in it (i.e. if it may simply code a sufficiently large Turing degree), arithmetic definability corresponds to membership in the minimal transitive model of ZFC\P + "$\omega_k$ exists" for some $k < \omega$; this is related to the complexity of Borel determinacy.

**Axiomatization**

While computability using pattern enumerators is $\Sigma_1^2$ definable (with Definition 4.2), it is insightful to axiomatize its key properties instead, as follows. We have not determined the strength or consistency of the axioms; we do not have a proof that bounded shuffles (below) permit the required erasure. We first axiomatize $g$ analogously to axiomatization of a sufficiently long fast-growing sequence in Section 4 (not included in the axioms below), with $g$ being sufficiently good relative to the enumerator or list of enumerators $f$ (unshuffled) (except that for the finite version, the cut-off point will be above $\max(g)$ (and thus could code $g$)). The below only assumes $g(n) \gg n$ and that $g$ is monotonic, with sharper bounds possible if $g(n+1) \gg g(n)$. The axioms are valid only for the domain of $g$. A level 1 infinite sequence appears to have the minimum expressiveness required (including for the transfinite extension below). However, the axioms below appear to depend on a level 2 sequence or higher to ensure that well-foundness as tested using $g$ agrees with well-foundness as tested using pattern enumerator(s). For a level 1 infinite sequence, we would also want that for every $s$ there is $t < g(s)$ such that $f(t)$ is the hyperjump of $f(s)$ as tested using $g$ (this also applies to the transfinite extension below). Below, $g$ as a numeric function refers to the lowest level. Also note that the section 4 axiomatization is for the nonlocal version.

As before, a *bounded shift* of $g$ above $p$ is any $g'$ with $g(i+p) \leq g'(i) \leq g(g(i+p))$ (for multilevel $g$, each level is shifted separately).

A *bounded shuffle* above $p$ (applied to the indices of $f$) is any permutation $Q$ with (for every $n$) $Q(n) < g(n+p)$ and $Q^{-1}(n) < g(n+p)$.

A computation $P$ using $p$ and $f(q)$ (both unaltered by the shift/shuffle) is *invariant* up to $n$ iff for every $m < n$, $P(m)$ is the same for every bounded shift and shuffle above $\max(p,q)$ and without the computation reaching the cut-off (if used) of the enumerator..

*Axioms:*
<u>Join axiom:</u> $f(i)$ join $f(j)$ is $f(k)$ for some $k < g(i+j)$.
<u>Completeness axiom/schema:</u> If a computation $P$ (using $p$ and $f(q)$) is invariant up to $n$ then there is $t < g(p+q)$ with $f(t)$ agreeing with $P$ up to $n$.
*Notes:*
* Since the tree of possibilities is finitely branching, the invariance for a particular example is recognizable.
* $g$ is accessible as a predicate, so a computation cannot test whether the end of $g$ (for the finite version) was reached.
* The completeness axiom is a schema over primitive recursive (or other desired complexity) $P$. Alternatively (if we want there to be a particular example), we can use a quantifier over $P$ (and $t < g(p+q+\text{code}(P))$), but without using a schema, the provable degree of completeness will be arbitrary for $g(i)$ for small $i$ (and for the finite version of the enumerator, arbitrary about how far above $g'$ we can go).

The completeness axiom states that if the computation is known (in a certain manner) to be invariant of the particular choice of $f$ and $g$, then the output is already a pattern. We do not have a predicate for pattern enumerators, but the expectation is that the shifts and shuffles can erase the extra information from the choice of $f$ and $g$. The exact treatment of 'bounded' (above) appears unimportant as long as it allows the erasure while keeping $g'$ (or some derivative) sufficiently long and fast-growing relative to the shuffled $f$.

**Transfinite Extension**

To go further, a sufficiently closed pattern enumerator of ordinal length will be a transfinite sequence of sufficiently closed pattern enumerators, each including a join of the previous ones as a pattern. Actually, the use of 'ordinal' (as opposed to just ordered) is superfluous since sufficient closure prevents an infinite descending sequence. It will be coded as $f(i)(j)(n)$ - character $n$ of pattern $j$ from enumerator $i$. The ordering of patterns and of enumerators is permitted to be arbitrary. For the ordinal length, one choice is to require the ordinal to be sufficiently closed relative to the notion of the pattern enumerator. Definition 4.2 works here (using the club filter for ordinals and using sets of reals as notions of pattern enumerators). We use $i \prec j$ for $f(i)$ being a pattern in $f(j)$. We use $g$ as before.

A bounded shuffle and invariance will be as before, except as follows:
* Each enumerator is shuffled independently. For enumerator $f(i)$, a shuffle above $p$ is above $\max(i, p)$.
* Require the transfinite join axiom (below) to hold after the shuffle.
* Also apply a bounded shuffle to the indices of the enumerators.
* An enumerator used as a parameter will also be shuffled.

*Axioms:*
* <u>Join:</u> $f(i)(j)$ join $f(i)(j')$ is $f(i)(k)$ for some $k < g(i + j + j')$.
* <u>Transfinite join:</u> Each enumerator $f(i)$ includes a pattern enumerating all other enumerators that do not include $f(i)$ as a pattern.
* <u>Completeness:</u> If a computation $P$ (using $p$, $f(i)(q)$ and $f(j)$) is invariant up to $n$, then there is $t < g(i + j + p + q)$ with $f(i)(t)$ agreeing with $P$ up to $n$.

*Notes:*
* See the note above about formalizing $P$ (including use of $\text{code}(P)$).
* Well-foundness of '$\prec$' should be provable (if using a level 2 $g$ or higher). If using a level 1 $g$, we can add a well-foundness axiom: $\exists n \forall x_1, \ldots, x_n \, (x_i < g^i(0)) \, \exists j \, x_{j+1} \not\prec x_j$. It is unclear if stronger axioms are needed.
* Using $g$, transfinite join can be a $\Sigma_1^0$ (or even bounded by $g^5(0)$) statement instead of a $\Pi_5^0$ statement.

Here is one choice of the ordinal closure axiom (different choices give different strengths):
There is an index $i < g(g(0))$ with $\forall i' < g(0) \, i' \prec i$ such that for every computation $P$ (using $p$, $f(j)$ and $f(k)(q)$ with $k \prec i$) that is invariant up to $n$, there is $l \prec i$ and $t$ with $l + t < g(i + j + p + k + \text{code}(P))$ with $f(l)(t)$ agreeing with $P$ up to $n$.
*Notes:*
* The condition on $i$ (before 'such that') means that the ordinal corresponding to $i$ can be arbitrarily large (we assume that the schema for $0 \ll g(0)$ is not cut off before it reaches the condition of the axiom). One can can create a stronger condition on $i$ using a schema that invariant properties reflect onto $i$ ($i$ dependent on the schema instance).
* This axiom should correspond to existence of an arbitrarily large ordinal $\alpha$ (represented by $i$ and thus within the range of '$\prec$') that becomes a cardinal in the model after higher level pattern enumerators are replaced with their invariant properties.
* We can go further by using indiscernible ordinals. A weak version is to strengthen the above axiom by having a tuple (using a schema or a quantifier over tuple length,

with the tuple length included in the bounds) of enumerators such that if a computation (using a subtuple as a parameter) is for every subtuple invariant (up to $n$), then the pattern $f(l)(t)$ (above) is independent (up to $n$) of the subtuple, other than the length and '$\prec$'-ordering of the subtuple, provided that $k \prec i \wedge j \prec i$ for every index $i$ in the subtuple (the condition on $j$ is not necessary for $l \prec i$).
\* One can extend the language further by including a sufficiently small club of countable ordinals (or a fast-growing tuple or sequence of such clubs) as a predicate and basing indices on a generic well-ordering of $\omega_1$, but it is unclear whether we have enough absoluteness to get a natural theory.

If consistent, these axioms appear to capture the basic form of the system, and perhaps the large cardinal strength (corresponding to enough Woodin cardinals). What remains is proving the strength of the system and finding natural axioms that capture appropriate large cardinal strength and make the theory reasonably complete for the results of invariant computations.

### 7.3 Using Countable Infinity

We can also consider what can be expressed using real numbers but without using an uncountable infinity. One way to increase expressiveness is to allow an infinite alternation of real quantifiers. Formally, $\exists X_1 \forall X_2 \ldots$ can be defined as true when the first player has a winning strategy in the game where the players takes turns to play $X_i$ and the goal of the first player is to satisfy the formula after the quantifiers. Under reasonable axioms, this leads to expressiveness beyond L($\mathbb{R}$). For every countable ordinal $\alpha$ (at least if $\alpha$ is not too large), real numbers definable from a countable ordinal using $1 + \alpha + 1$ real quantifiers are precisely those that are in the minimal iterable inner model with $\alpha$ Woodin cardinals. Just $\omega + 2$ real quantifiers go beyond first order definability in L($\mathbb{R}$).

This definition uses strategies, which are formally uncountable. However, we can avoid the uncountable by introducing a notion of $Y$ having a sufficiently large complexity (under Turing reducibility) relative to $X$. As before, vagueness is avoided because every sufficiently strict notion works here. We obtain the necessary complexity by using a sequence of $\alpha$ real numbers, each sufficiently complex relative to a join of the previous real numbers. Statements with infinite alternations of quantifiers can be decided by restricting the complexity of the strategies: For $\gamma$ levels of play, the strategies must be recursive in the $\gamma$'s member of the sequence. By determinacy and by sufficiency of the sequence, this does not affect who wins. We can go even further by using sequences with ordinal length sufficiently long relative to themselves. Reaching higher complexity is an endless endeavor.